\documentclass
[11pt]{article}
\usepackage{tipa}
\usepackage{amssymb}
\usepackage{mathrsfs}
\usepackage{amsfonts}

\usepackage{graphicx}

\include{monontone_ref}

\makeindex             

\def\go{\rightarrow}

\begin{document}
\title{\textbf{Asymptotic Behavior of Coarse-grained Models for Opinion
Dynamics on Large Networks}}
\author{Chjan C. Lim and Weituo Zhang \\
Mathematical Sciences, RPI, Troy, NY 12180, USA\\
email: limc@rpi.edu} \maketitle
\begin{abstract}
In this paper, we propose a general mathematical framework to
represent many multi-agent signalling systems in recent works. Our
goal is to apply previous results in monotonicity to this class of
systems and study their asymptotic behavior. Hence we introduce a
suitable partial order for these systems and prove nontrivial
extensions of previous results on monotonicity. We also derive a
convenient sufficient condition for a signalling system to be
monotone and test our condition on the Naming Games, NG and K-NG on
complete networks both with and without committed agents. We also
give a counter example which fails to satisfy our condition. Next we
further extend our conclusions to systems on sparse random networks.
Finally we discuss several meaningful consequences of monotonicity
which narrows down the possible asymptotic behavior of signalling
systems in mathematical sociology and network science.
\end{abstract}

\section{Introduction}

We introduce a mathematical framework comprehensive enough to
represent many of the multi-agent signalling networks in recent
works such as the generalized Naming Games (NG) \cite{Xie11,
ZhangLim, Yosef, Baronchelli12} and the Voter models used in
modeling of ad-hoc wireless networks, agreement on a small list of
tags in the WWW and social opinion dynamics in forums \cite{Lim13}
and elections.

Our aim is to show that the asymptotic behavior of a proper subset
of these signalling games on very large networks (where the effects
of demographic noise can be safely neglected), can be described
rigorously as a system of nonlinear coupled ODEs with the
Monotonicity Property (MP) in phase-space consisting of the
coarse-grained population fractions of each type. The nonempty
complement of the monotone signalling systems include two main
types, namely, (1) the signalling systems where the associated
random walk models are Martingales, that is, the deterministic
drifts are zero everywhere in phase space; in such systems, which
includes all the Voter models, the asymptotic behavior is diffusion
- driven \cite{Star}, \cite{Lim13}, and (II) signalling systems
which have nonzero drift almost everywhere in phase space but whose
nonlinear ODEs are not monotone; we will give such a counterexample
below.

Previous works on Monotonicity include the works of \cite{Smith,
Smith2} which are mainly used in applications to population genetics
and mathematical ecology \cite{Levin}. In this paper, our main aim
is to extend the applications of monotonicity to information sharing
and opinion dynamics on large social networks of interacting
(signalling) agents consisting of a finite, but possibly large,
number of opinion or information types. The mathematical framework
given in detail below includes all known multi-agent
binary-signalling networks, that is, where the signals are
restricted to the binary symbols ${A,B}$. There are many ways to see
why the binary-signals case represents the end-game in the final
stages of signalling networks based on more than two symbols.
Concrete evidence of this has been discussed in \cite{Zhangthesis,
Lim13} where expected times to multi-consensus for a symmetric
multi-opinions Voter model was calculated in terms of the bottleneck
times when exactly two opinions remain viable in the game
\cite{Lim13}.

In the process of working out the mathematical and socio-physical
consequences \cite{Granovetter} \cite{Schelling1978} of monotonicity
\cite{Smith}, we will employ new versions of partial order that are
explicitly designed for signalling networks. In particular, we focus
on a new way to overcome a previous obstacle to the direct
applications of Monotonicity results, namely, a mismatch in the
cardinality of the set of local (edge or node) spin types and the
number of independent (edge or node) population types. This mismatch
arises specifically in attempts to apply the traditional
Monotonicity results to signalling systems on large social
(especially random) networks other than the complete graph. Along
the way, we will also prove nontrivial extensions of the original
monotonicity theorems in \cite{Smith, Smith2}.

\section{Framework-General signalling system and its macrostate space}

Consider a signalling network consisting of N agents, each of which
is assigned with a spin $s_i$ taking value from the spin state space
$\Gamma=\{\gamma_1,\gamma_2,...,\gamma_K\}$. The (micro)state of the
system is fully described by the vector of spins
$\vec{S}=(s_1,...,s_i,...,s_N)$.

In each time step, a speaker and a listener are randomly selected.
The speaker sends a message containing 1 bit of information to the
listener, and the listener changes its state according to the
message. In this context, the message can be ``A'' or ``B''. The
probability for the selected speaker in state $\gamma_k$ sending a
message ``A'' is $\alpha_k$. Let
$\vec{\alpha}=(\alpha_1,...,\alpha_K)^T$ ($\alpha_k\in[0,1]$).

When listener receive a message A (resp. B), the transition matrix
of the listener's state is $G_A$ (resp. $G_B$), the $i-j$th entry of
which is $P(s_i\go s_j|A)$ for $G_A$.
Both $G_A$ and $G_B$ are constant matrices and the signalling
dynamics of the network is fully governed by $\vec{\alpha}$, $G_A$
and $G_B$.

A natural macrostate representation of this system is given by
$\vec{n}=(n_1,...,n_k,...,n_K)^T$ where $n_k$ is the fraction of nodes in spin state
$\gamma_k$ ($k=1,...,K$). The space of all possible macrostates is denoted by $M$, which is a simplex:
\begin{equation}
\left\{\begin{array}{l}
       \sum_{i=1}^d n_i=1\\
       n_i \geq 0 \ \ (i=1,...,d)
       \end{array}\right.
\end{equation}

We refer to the vertices of $M$ as pure macrostates in which all
nodes of the system stay in the same spin state. For later use, let
$\sigma:\Gamma\go M$ map a spin state to its corresponding pure
macrostate. We rewrite the macrostate as a linear combination of
pure macrostates.
$$\vec{n}=(n_1,...,n_k,...,n_K)=\sum_{k=1}^K n_k\sigma(\gamma_k)$$

\section{Partial order}
We begin with $\{\Gamma,\prec\}$, the spin state space ordered by a partial order relation ``$\prec$'', satisfying:
\begin{equation}
\left\{\begin{array}{l}
      \mbox{reflexivity}: \gamma\prec\gamma \\
      \mbox{symmetry}: (\gamma\prec\gamma') \wedge (\gamma'\prec\gamma) \Rightarrow \gamma=\gamma' (are equal) \\
      \mbox{transitivity}:(\gamma_1\prec\gamma_2) \wedge (\gamma_2\prec\gamma_3) \Rightarrow \gamma_1\prec\gamma_3
       \end{array}\right.
\end{equation}
In later sections we will discuss what further conditions this
partial order should satisfy. We induce a partial order relation
between the pure macrostates from $\{\Gamma,\prec\}$ through
$\sigma:\Gamma\go M$:
$$\sigma(\gamma)\prec \sigma(\gamma') \iff   \gamma\prec\gamma'$$
We use the same notation for the partial orders in both spaces with no ambiguity.
Next we extend the partial order over $\Gamma$ through affine combination. We require
$$\vec{n_1}\prec \vec{n_1}', \vec{n_2}\prec \vec{n_2}' \Rightarrow  L(\vec{n_1},\vec{n_2})\prec L(\vec{n_1}',\vec{n_2}').$$

Here $L(\vec{n_1},\vec{n_2})=a \vec{n_1}+ b \vec{n_2}$ is an
affine combination of $\vec{n_1}, \vec{n_2}$ which satisfies $0\leq a,b \leq 1$
and $a+b=1$ so that $L(\vec{n_1},\vec{n_2})\in M$ whenever $\vec{n_1},\vec{n_2}\in M$.
Now the extended partial order, $\vec{n}\prec\vec{n}'$ if and only if $\vec{n}'-\vec{n}$ can be represented as
$$\vec{n}'-\vec{n}=\sum_i \lambda_i (\sigma(\gamma_i')-\sigma(\gamma_i))$$
The above summation is over all possible ordered pure macrostate pairs $\sigma(\gamma_i)\prec \sigma(\gamma_i')$
and $\lambda_i\geq 0$ for all $i$. Note that some terms in the summation may be further decomposed into
combination of other ordered pairs, so we consider the independent set
$$B=\{\sigma(\gamma_i')-\sigma(\gamma_i)| \gamma_i\prec \gamma_i'  \mbox{, and } \nexists \gamma'' st. \  \gamma_i\prec \gamma'' \prec\gamma_i'\}.$$

Let $TM$ be the tangent space of macrostate space $M$. Then
$C_+\subset TM$ is the cone determined by the non-negatively linear
combination of set $B$. The definition of partial order here can be
described as

\begin{equation}
\vec{n}\prec \vec{n}' \Longleftrightarrow  \vec{n}'-\vec{n}\in C_+
\end{equation}

If $|B|$, the cardinality of $B$, equals $dim(TM)$, the dimension of
$TM$, the partial order we described here can be considered as a
standard type K partial order in Cartesian coordinate system
($x\prec y \iff y-x\in \mathbb{R}_+^n$) after suitable linear
transformation \cite{Smith}. However, it is possible that $|B|\neq
dim(TM)$ regarding some specific ``$\prec$'' on $\Gamma$ (we will
provide some examples later), therefore our partial order given by
$C_+$ is a more general definition than the standard ones in the
literature. It is clear that this more general form is needed for
the applications in this paper to some well-known network games
\cite{Zhangthesis}.

\section{Macrostate dynamics of signalling system on complete networks}
On a complete network, assume at time step $t$, the macrostate is
$\vec{n}(t)$, then at time step $t+1$, the expected change of
macrostate is
$$E[\vec{n}(t+1)-\vec{n}(t)]={1\over N}\left[pG_A+(1-p)G_B-I\right] \vec{n}(t),$$
where p is the overall probability for a message to be ``A'', given
by
$$p=\vec{\alpha}^T\vec{n}.$$

With standard time scaling $dt=1/N$, we obtain the mean field equation for the evolution of signalling system:
$${d\vec{n}\over dt}=f(\vec{n})=Q(\vec{n})\vec{n}(t)=\left[pG_A+(1-p)G_B-I\right] \vec{n}(t).$$
Here $f$ maps $M$ to its tangent space $TM$. Denote $\phi_t$ as the semiflow
which give the solution of the ODE, $\vec{n}(t)=\phi_t(\vec{n}_0)$.

A system is said to be order-preserving or monotone,
if for $\forall \vec{n}\prec \vec{n}'$ and $\forall t>0$, we have $\phi_t(\vec{n})\prec \phi_t(\vec{n}')$.

We now give a condition for monotonicity, which is an analogue of
the Kamke condition \cite{Smith}. Assume
$B=\{\vec{e}_1,...,\vec{e}_{|B|}\}$. $f$ is said to be
\textbf{{\emph{type C}}} if for each $k\in\{1,...,|B|\}$, for
$\forall \vec{n}\prec \vec{n}'$ such that
$\vec{n}'-\vec{n}=\sum_{i\neq k}a_i \vec{e}_i$ ($a_i\geq 0$), there
exists a representation of $f(\vec{n}')-f(\vec{n})=
\sum_{i=1}^{|B|}b_i \vec{e}_i$ in which $b_k\geq 0$. Note that when
$|B|$ is greater than $dim(TM)$, the representation here may not be
unique, but when $|B|<dim(TM)$, the representations may not exist
regardless the sign of the coefficients.

{\bf Proposition 1}(type C condition): The system ${d\vec{n}\over dt}=f(\vec{n})$ is monotone if and only if it is type C.

The proof is straight forward by contradiction and very similar to
that of the Kamke condition. If the type C condition does not hold,
then by continuity of $f$, there exists an $\epsilon>0$ such that
$\phi_{\epsilon}(\vec{n})\succ \phi_{\epsilon}(\vec{n}')$ which
violates the monotone property.

As the Kamke condition \cite{Kamke,Smith} can be expressed in terms of
partial derivatives, the type C condition has an expression in terms
of directional derivatives.

{\bf Proposition 2}: The system ${d\vec{n}\over dt}=f(\vec{n})$ is monotone if and only if (A): for $\forall \vec{n}\in Int(M)$ and $\forall k\in\{1,...,|B|\}$, there exists a representation ${d\over d\epsilon} f(\vec{n}+\vec{e_k})=\sum_{i=1}^{|B|}b_i \vec{e_i}$ s.t. for $\forall i\neq k$, $b_i\geq 0$.

{\bf Proof}: For $\forall \vec{n}\prec \vec{n}'$ such that
$\vec{n}'-\vec{n}=\sum_{i\neq k}a_i \vec{e}_i$ ($a_i\geq 0$) and the
representation of $f(\vec{n}')-f(\vec{n})=
\sum_{i=1}^{|B|}b_i(\vec{n}) \vec{e}_i$ holds with the sign of $b_i$
undecided,

\begin{eqnarray}
f(\vec{n}')-f(\vec{n})&=&\int_0^1 {d\over d\lambda}\vec{f}(\vec{n}+\lambda \sum_{i\neq k}a_i \vec{e}_i)d\lambda \label{eqn1}\\
                     &=&\sum_{j\neq k} \int_0^1 a_i {d\over d\epsilon}\vec{f}(\vec{n}+\lambda \sum_{i\neq k}a_i \vec{e}_i+ \epsilon \vec{e}_j)d\lambda
\end{eqnarray}

(a) (A)$\Rightarrow$type C: If conditon (A) holds, then ${d\over
d\epsilon}\vec{f}(\vec{n}+\lambda \sum_{i\neq k}a_i \vec{e}_i+
\epsilon \vec{e}_j)=\sum_{i=1}^{|B|}b_i^{(j)}(\lambda) \vec{e_i}$
such that for $\forall i\neq j$, $b_i\geq 0$. We point-wisely
replace the directional derivative by its linear representation,
therefore
\begin{equation}
f(\vec{n}')-f(\vec{n})=\sum_{j\neq k} \int_0^1 a_i \sum_{i=1}^{|B|}b_i^{(j)}(\lambda) \vec{e_i} d\lambda. \label{eqn2}
\end{equation}
Summing all the coefficients before $\vec{e}_k$ we get $\sum_{j\neq
k} a_k \int_0^1 b_k^{(j)}(\lambda)d\lambda \geq 0$, so $f$ is type
C.

(b) Type C$\Rightarrow$(A): Suppose for some $\vec{n}_0\in Int(M)$,
${d\over d\epsilon}
f(\vec{n}_0+\epsilon\vec{e_j})=\sum_{i=1}^{|B|}b_i \vec{e_i}$,
$b_k<0$ ($k\neq j$). By continuity of $f$, there exists a small
enough $\delta_0>0$ such that $b_k<0$ also holds for
$\vec{n}=\vec{n}_0+\delta \vec{e_j}$ when $0<\delta\leq \delta_0$.
Applying Eq.~(\ref{eqn2}) $f(\vec{n}')-f(\vec{n})=\delta_0
\sum_{i=1}^{|B|}\left(\int_0^1 b_i^{(j)}(\lambda) d\lambda\right)
\vec{e_i}$. Since $b_k^{(j)}(\lambda)<0$ for $0\leq\lambda\leq 1$,
the coefficient before $\vec{e_k}$ is $\delta_0
\sum_{i=1}^{|B|}\left(\int_0^1 b_k^{(j)}(\lambda) d\lambda\right)<0$
contradicting the type C condition. QED

For the signalling system on a complete network, we have
$f(\vec{n})=\left[pG_A+(1-p)G_B-I\right] \vec{n}(t)$, thus,
\begin{eqnarray*}
{d\over d\epsilon} f(\vec{n}+\epsilon \vec{e}_k)&=& (pG_A+(1-p)G_B-I)\vec{e}_k+ \left[{dp(\vec{n}+\epsilon\vec{e}_k)\over d\epsilon}{d\over dp}(pG_A+(1-p)G_B)\vec{n} \right] \\
&=& pG_A\vec{e}_k+(1-p)G_B\vec{e}_k-\vec{e}_k+(\vec{\alpha}^T\vec{e}_k)(G_A-G_B)\vec{n}\\
\end{eqnarray*}

Now we derive a sufficient conditions for condition (A). According to Proposition 2, only $b_i$ $(i\neq k)$ affects the monotonicity, so the third term in the last expression ($-\vec{e}_k$) can be neglected. If we require the other three terms be positive in terms of $\prec$, we obtain the following theorem.

\textbf{Theorem 1:} A signalling dynamics on complete network
governed by $\vec{\alpha}$, $G_A$, $G_B$ is monotone if it satisfies
\begin{equation}
\left\{ \begin{array}{l}
       (a).  \forall \gamma\in \Gamma, \  G_B\sigma(\gamma)\prec G_A\sigma(\gamma)\\
       (b).  \gamma\prec \gamma'  \Rightarrow \vec{\alpha}^T\sigma(\gamma)<\vec{\alpha}^T\sigma(\gamma')\\
       (c). \gamma\prec \gamma' \Rightarrow G_A \sigma(\gamma)\prec G_A \sigma(\gamma'), G_B \sigma(\gamma)\prec G_B \sigma(\gamma')
       \end{array}\right.
\end{equation}

This theorem gives a convenient sufficient condition on the
essential components of the signalling system in order for it to be
monotone. In this theorem, (a) fixes the preferred (greater) one
between A and B to orient the partial order, that is, the listener
will switch to a greater state according to this partial order when
receiving A than when receiving B. Condition (b) says the speaker in
a greater state has more probability to send a message A. One can
switch the roles of A, B in the conditions (a), (b). Condition (c)
says $G_A$ and $G_B$ preserve the partial order. In other words, we
can now determine directly whether a signalling system is monotone
according to the properties of its three governing elements
$\vec{\alpha}, G_A, G_B$.

For particular applications of these sufficient conditions to social
networks of signalling agents, we can imbue the partial order with a
moral value system or a utility function \cite{Granovetter}

\section{Examples of binary signalling system on complete graph}
\subsection{Binary Listener-only Naming Game (LO-NG) with and without committed agents}
In this case \cite{Xie11}, \cite{ZhangLim}, the spin states
$\Gamma=\{A, AB, B\}$, and a macrostate is given by the
corresponding populations $\vec{n}=(n_A,n_{AB},n_B)^T$. Here, the
governing elements are given by the following vector which fixes the
probabilities of sending the symbol $A$ when in the the associated
sin states and a pair of transition matrices, which define the
transition probabilities upon receiving a $A$ (resp. $B$) symbol:
$$\vec{\alpha}=(1,1/2,0)^T,$$
$$G_A=\left( \begin{array}{ccc}
          1 & 1 & 0\\
          0 & 0 & 1\\
          0 & 0 & 0
        \end{array}\right), G_B=\left( \begin{array}{ccc}
          0 & 0 & 0\\
          1 & 0 & 0\\
          0 & 1 & 1
        \end{array}\right). $$
Thus, the governing elements of a signalling system of networked
agents have two natural parts, namely,  (I) the probabilities for
sending a symbol, and (II) the transition probabilities to the next
local state (opinion) on receiving a symbol.

After choosing the opinion $A$ to be "morally superior" say, and
inducing from the implicit order in $\vec{\alpha}$, the suitable
partial order over $\Gamma$ that satisfies condition (b) in the main
theorem 1, should be $B\prec AB \prec A$ which . Therefore, the
independent subset of pairs, identified in the general mathematical
framework given above, is $B=\{\sigma(A)-\sigma(AB),
\sigma(AB)-\sigma(B) \}=\{(1,-1,0),(0,1,-1)\}$, where the second
equality should obviously mean that the negative ones denote the
second element / term of the differences in the independent set $B$,
that is, the last set of 3-vectors are incidence vectors for the
differences in $B$ viewed as directed edges. Thus, this partial
order can be represented by a directed graph, in which $\Gamma$
gives the vertices and $B$ gives the directed links, as in the
figure below.

\begin{figure}[!hbtp]
\begin{center}
\caption{Partial order of NG}
\vspace{0cm}
  \includegraphics[width=0.4\textwidth]{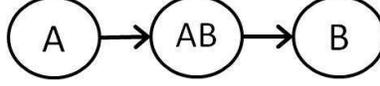}\\
 \end{center}
\end{figure}

Next we check the conditions (a) and (c) in Theorem 1. In
Table~\ref{tab_NG}, comparing two rows, we verify (a)
$G_B(\sigma(X))\leq G_A(\sigma(X))$ ($X\in\{A, AB, B\}$); comparing
three columns, we verify (c) $G_X(\sigma(B))\prec
G_X(\sigma(AB))\prec G_X(\sigma(A))$ ($X\in\{A,B\}$).

\begin{table}[!h]
\begin{center}
\caption{$G_A$, $G_B$ of NG}
\label{tab_NG}
\begin{tabular}{cccc}\hline
                   & A & AB & B\\   \hline
$G_A(\sigma(\cdot))$ & $\sigma(A)$ & $\sigma(A)$ & $\sigma(AB)$\\
$G_B(\sigma(\cdot))$ & $\sigma(AB)$  & $\sigma(B)$ & $\sigma(B)$\\
\hline
       \end{tabular}
       \end{center}
       \end{table}

In the case with committed agents both on A and B, $\Gamma=\{C_A, A,
AB, B, C_B\}$; $C_A$ and $C_B$ denote the spin state committed in A
and in B respectively and the populations macrostate is given by
$\vec{n}=(n_{C_A}, n_A,n_{AB},n_B, n_{C_B})^T$. The governing
elements in this case are
$$\vec{\alpha_c}=(1, 1, 1/2, 0, 0)^T,$$
$$G_{A,c}=\left( \begin{array}{ccc}
        1&  &\\
          &G_A&\\
          && 1\\
        \end{array}\right), G_{B,c}=\left( \begin{array}{ccc}
          1& &\\
           &G_B&\\
            & & 1\\
        \end{array}\right). $$
The suitable partial order is shown in the directed graph below
(just add two disconnected points to that of the non-committed or
symmetric LO-NG case). Committed agents never changes its state, so
$G_A(\sigma{X})=G_B(\sigma{X})=\sigma{X}$ ($X\in \{C_A, C_B\},$)
which therefore satisfy condition (a) in Theorem 1. Conditions (b)
and (c) only affect the pairs ordered by the partial order; since
the committed case does not introduce any new ordered pairs w.r.t.
the non-committed case, and $\vec{\alpha_c}, G_{A,c}, G_{B,c}$
restricted to the non-committed pure macrostates are exactly
$\vec{\alpha}, G_{A}, G_{B}$, hence conditions (b) and (c) n Theorem
1 are satisfied.

From here, we can easily get the following corollary which has
applications to the scenario in \cite{korniss} and \cite{swami}:

{\bf Corollary:} If a signalling system without committed agents is
monotone, adding committed agents into this signalling system does
not change the monotonicity.

\begin{figure}[!hbtp]
\begin{center}
\caption{Partial order of NG with committed agents}
\vspace{0cm}
  \includegraphics[width=0.5\textwidth]{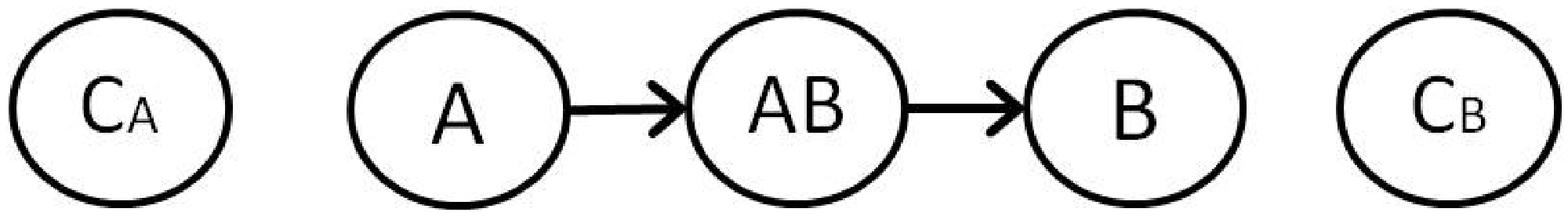}\\
 \end{center}
\end{figure}

\subsection{K-NG}
\begin{figure}[!hbtp]
\begin{center}
\caption{Partial order of K-NG}
\label{fig_KNG}
\vspace{0cm}
  \includegraphics[width=0.5\textwidth]{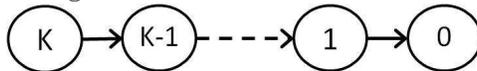}\\
 \end{center}
\end{figure}

In the one-parameter family of listener-only Naming Games, K-NG
\cite{Yosef}, where $K$ presents the stubbornness of agents to full
conversion from $B$ to the $A$ opinion, there are $K+1$ spin states
$\Gamma=\{0, 1,..., K\}.$  As in the original LO-NG case which
corresponds to the value $K=2$ in the $K-NG$ family of models, we
firstly find a suitable partial order satisfying condition (b) in
Theorem 1, that agrees with the given probability vector for sending
the $A$ symbol, $\vec{\alpha}=(0, 1/K, ..., i/K, ..., 1)^T$. This
partial order is shown in the Fig.~\ref{fig_KNG}. Next we find the
independent subset $B=\{\sigma(k+1)-\sigma(k)| k=0,...,K-1\}$.

Since the $K-NG$ family of NG models have the following property:
$$\sigma(max(k-1,0))\prec G_B(\sigma(k))\prec \sigma(k)\prec G_A(\sigma(k))\prec \sigma(max(k+1,K)),$$
conditions (a) and (c) in Theorem 1 follow easily.

\subsection{Counter-example which is not monotone}
Next we provide an explicit example of a binary signalling system
which is not monotone, thus establishing that the above definitions
in our mathematical framework has non-vacuous complement. This
example of a signalling system on two symbols has the same spin
state space $\Gamma$ and macrostate representation as the LO-NG
model discussed previously. However, its governing elements (the
probabilities for sending $A$ and the spin transition probabilities
upon receiving the symbols $A$ or $B$) differ from the LO-NG in the
matrix $G_B$:
$$\vec{\alpha}=(1,1/2,0)^T,$$
$$G_A=\left( \begin{array}{ccc}
          1 & 1 & 0\\
          0 & 0 & 1\\
          0 & 0 & 0
        \end{array}\right), G_B=\left( \begin{array}{ccc}
          0 & 0 & 0\\
          0 & 1 & 0\\
          1 & 0 & 1
        \end{array}\right). $$

We claim there does not exist a non-trivial partial order in this
example, where by nontrivial we mean that the partial order has a
nonempty independent subset $B\neq\emptyset.$ The proof of this
claim follows:

{\bf Proof}: According to condition (b) in Theorem 1, there are at
most three possible elements in $B$ for this example, i.e. $B\subset
\{\vec{e}_1=\sigma(A)-\sigma(B), \vec{e}_2=\sigma(AB)-\sigma(B),
\vec{e}_3=\sigma(A)-\sigma(AB)\} $. Since $G_A \vec{e}_1=\vec{e}_3,
G_A \vec{e}_2=\vec{e}_3$, by (c) in Theorem 1, $\vec{e}_1\in
B\Rightarrow \vec{e}_3\in B$ and $\vec{e}_2\in B\Rightarrow
\vec{e}_3\in B$. If $B\neq\emptyset$, then $\vec{e}_3\in B$.
However, $G_B \vec{e}_3=\sigma(B)-\sigma(AB)$, by (c) again,
$\sigma(AB)\prec \sigma(B)$ which contradicts condition (b) in
Theorem 1. QED

The solution trajectories of this system mapped into 2D space $(n_A,
n_B)$ is shown below. This figure gives us a concrete idea of one
type of non-monotone signalling dynamics on two symbols.

\begin{figure}[!hbtp]
\begin{center}
\caption{Trajectories of counter example}
  \includegraphics[width=0.6\textwidth]{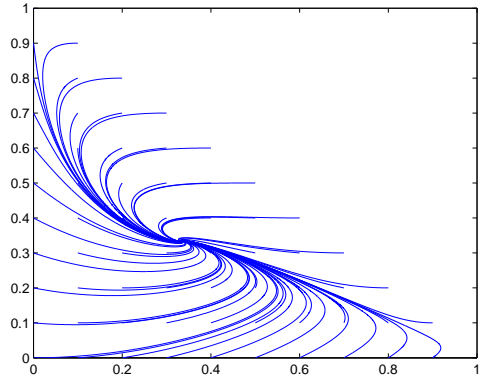}\\
 \end{center}
\end{figure}

\section{On sparse random networks}
Under a homogeneous pairwise mean field assumption, the
coarse-grained dynamics of the binary Naming Game on a random sparse
network with average degree $<k>$ is governed by the 6-dimensional
ODEs \cite{Zhangcomple12} \cite{ZhangLimBolek12}:
\begin{eqnarray}
{d\over dt}\vec{l}=2\left[{1\over <k>}D+({<k>-1\over <k>})R\right]\vec{l}.
                  \label{eq_NGER}
\end{eqnarray}
for macrostate
$\vec{l}=[l_{A-A},l_{A-B},l_{A-AB},l_{B-B},l_{B-AB},l_{AB-AB}]^T$ of
link-types population fractions. For a general signalling system of
the above type, we can obtain a similar ODE system using exactly the
same approach as in \cite{ZhangLimBolek12}. The type of a link in
this ODEs is ultimately given by the opinions or node-spins of its
two ends ($\gamma_i-\gamma_j$) regardless of the order but it is
more convenient to work with link-based macrostates; a
transformation between the link-based macrostate and the node-based
macrostate $\vec{n}$ is given in \cite{ZhangLimBolek12}.
Changes of $\vec{l}$ come from two parts: the \emph{\textbf{direct
change}} and the \emph{\textbf{related change}}. In each time step,
a realized listener-speaker pair of agents and therefore a link or
edge in the underlying random graph is selected. Then
\textbf{\emph{Direct change}} is the change of the selected link and
is given by $D\vec{l}$ where $D$ is the probability transition
matrix of the selected link - the $(i,j)$ entry of $D$, $D_{ij}$, is
the probability that a link of type j changes into type i in one
step given that the selected link is of type j. $D$ is a constant
matrix given by $\vec{\alpha}$, $G_A$, $G_B$ and the way of
selecting listeners and speakers.

Consider a link type $\gamma_1-\gamma_2$ and its corresponding pure
macrostate $\sigma(\gamma_1-\gamma_2)$ which, according to the
general mathematical framework in this paper, are now the basic
elements on which a partial order is defined. Then according to this
theory and the above 6-dimensional ODEs (with effectively a 5-dim
tangent space after reduction by one degrees of freedom),
\begin{eqnarray*}
D\sigma(\gamma_1-\gamma_2)&=& P(\gamma_1 \mbox{ is listener}) \left[p_1 \sigma(G_A \gamma_1-\gamma_2)+ (1-p_1)\sigma(G_B \gamma_1-\gamma_2)\right] \\
               &&+P(\gamma_2 \mbox{ is listener}) \left[p_2 \sigma(\gamma_1- G_A \gamma_2)+ (1-p_2)\sigma(\gamma_1-G_B \gamma_2)\right],
\end{eqnarray*}
where $p_1=\vec{\alpha}^T\sigma(\gamma_2)$,
$p_2=p_1=\vec{\alpha}^T\sigma(\gamma_1)$. The notations $\sigma(G_A
\gamma_1-\gamma_2)$ is defined by
$$\sigma(G_A \gamma_1- \gamma_2)=\sum_{i=1}^K P(G_A \sigma(\gamma_1)=\sigma(\gamma_i)) \sigma(\gamma_i-\gamma_2)$$

The related links are those incident at the listener other than the
selected link. Then the \emph{\textbf{Related change}} is the change
of the related links when the listener changes and is given by
$(<k>-1)R(\vec{l})\vec{l}$, where $<k>-1$ is the expected number of
related links. $R$ is a probability transition matrix that varies
according to the current macrostate. For Naming Game, $R$ is given
explicitly in \cite{ZhangLimBolek12}. For general signalling
systems, $R$ is given in detail later.

A natural partial order of link-based macrostates is induced from that of node-based macrostates: the link states
$$X-Y\prec X'-Y' \iff X\prec X' and Y\prec Y'.$$
The independent set of ordered pairs $B$ is thus
$$B=\{\sigma(X-Y')-\sigma(X-Y)| Y\prec Y' \}.$$
For example, the partial order and the set $B$ for NG on a sparse
random network is shown in the following figure. Note that it is an
example in which $|B|=6\geq dim(TM)=5$.
\begin{figure}[!hbtp]
\begin{center}
\caption{Partial order of NG on random sparse networks}
\vspace{2cm}
  \includegraphics[width=0.2\textwidth]{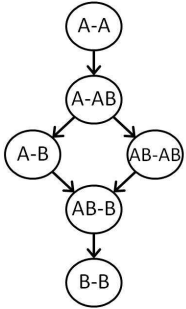}\\
 \end{center}
\end{figure}

The partial order in the link-based macrostate space is stronger
than that in the node-based macrostate space, i.e. if two node-based
macrostates are $\vec{n}, \vec{n}'$ and their respective linked
based macrostates $\vec{l}, \vec{l}'$, then $\vec{l}\prec
\vec{l}'\Rightarrow \vec{n}\prec \vec{n}'$ and the reverse does not
hold.

In Section 4, we proved a sufficient condition for a signalling
dynamics to be monotone on a complete network graph. We will show in
the following theorem, that the signalling dynamics satisfying this
condition will also be monotone on a random sparse network.

\textbf{Theorem 2:} If a signalling dynamics satisfies the conditions
in Theorem 1, then it is monotone on a random sparse network.

\textbf{Proof:} Firstly, we show that given the conditions in (6),
the direct change part $D\vec{l}$ satisfies condition (a). For each
$\vec{e}_k=\sigma(X-Y')-\sigma(X-Y)$,
\begin{eqnarray*}
D\vec{e}_k&=& D\sigma(X-Y')-D\sigma(X-Y)\\
          &=& P(X \mbox{ is listener}) [ \vec{\alpha}^T \sigma(Y)\left(\sigma(G_A X-Y')-\sigma(G_A X-Y)\right)\\
          &&+\left(1-\vec{\alpha}^T \sigma(Y)\right)\left(\sigma(G_B X-Y')-\sigma(G_B X-Y)\right)\\
          &&+ \left(\vec{\alpha}^T \sigma(Y')-\vec{\alpha}^T \sigma(Y)\right)\left(\sigma(G_A X-Y')-\sigma(G_B X-Y')\right)]\\
          &&+ P(Y \mbox{ or } Y' \mbox{ is listener}) [\vec{\alpha}^T \sigma(X)\left(\sigma(X-G_A Y')-\sigma(X-G_A Y)\right)\\
          &&+ \left(1-\vec{\alpha}^T \sigma(X)\right)\left(\sigma(X-G_B Y')-\sigma(X-G_B Y)\right)]
\end{eqnarray*}

By the definition of link-based partial order, $\sigma(G_A X-Y')-\sigma(G_A X-Y)\succ 0$, $\sigma(G_B X-Y')-\sigma(G_B X-Y)\succ 0$. According to condition (6), $\vec{\alpha}^T \sigma(Y')-\vec{\alpha}^T \sigma(Y)>0$, $\sigma(X-G_A Y')-\sigma(X-G_A Y)\succ 0$, $\sigma(X-G_B Y')-\sigma(X-G_B Y)\succ 0$. Therefore $D\vec{l}$ preserve the order of $l$

For the related change part, $R\vec{l}$, we represent $\vec{l}$ as a
suitably weighted symmetric adjacency matrix $M=M(\vec{l})$ labeled
by the node-based spin types,
\begin{equation}
\left\{\begin{array}{l}
      M_{ii}=l_{\gamma_i-\gamma_i},\\
      M_{ij}={1\over 2} l_{\gamma_i-\gamma_j} (i\neq j).
       \end{array}\right.
\end{equation}
where the row sum and column sum of $M(\vec{l})$ are the node-based
macrostate, $\vec{n}$. It is obvious that, $M(\vec{l})$ is a $1-1$
presentation for $\vec{l}$. $R(\vec{l})$ is given by the following
$$M(R(\vec{l})\vec{l'})=W(\vec{l}) M(\vec{l'}) W(\vec{l})^T,$$

where $W(\vec{l})$ is the transition matrix of spin states, i.e. the
entry of $W(\vec{l})$, $W_{ij}$ is the probability that a node in
spin state $\gamma_j$ changes into $\gamma_i$ given the link-based
macrostate $\vec{l}=(l_1,...,l_k,...,l_K)$,
$$W_{ij}=\sum_{k=1}^K P(\gamma_j\go \gamma_i| \mbox{ link type } \ k) l_k .$$

Then we prove the following lemma.

\textbf{Lemma 1:} There exists a unique decomposition of
$M(\vec{l})$, $M(\vec{l})=u+v  $. $u\in U=\{u|u=u^T, \exists
\mbox{column vector} \vec{m} \mbox{ s.t. } u=\vec{m}\otimes
\vec{m}^T\}$, $\otimes$ is the kronecker product . $v\in
V=\{v|v=v^T,  v\vec{\bold{1}}=\vec{0}, \vec{\bold{1}}^T
v=\vec{0}\}$, where $\vec{\bold{1}}$ is a column vector with
all $1$ entries.

\textbf{Proof:} Taking $\vec{m}=\vec{n}$ the node-based population
fractions, $u=\vec{n}\otimes \vec{n}^T$ is symmetric by
construction. The row sum and column sum of $u$ are both $\vec{n}$,
the same as $M(\vec{l})$. Therefore the row sum and column sum of
$v$ are $0$. As $M(\vec{l})$ and $u$ are symmetric, so is $v$.

\textbf{Lemma 2:} $U$ and $V$ are invariant space of the operator $R(\vec{l})$ for any macrostate $\vec{l}$.

\textbf{Proof:}
$\forall u=\vec{m}^T\otimes \vec{m}\in U$, $WuW^T=(W\vec{m})\otimes (W\vec{m})^T \in U$.
$\forall v\in V$, $v\vec{\bold{1}}=\vec{0}$, since $W^T\vec{\bold{1}}=\vec{\bold{1}}$, $WvW^T\vec{\bold{1}}=\vec{\bold{1}}^TWvW^T=\vec{0}$, therefore $WvW^T\in V$.

According to Lemma 1 and 2, $WM(\vec{l})W^T$ restricted on $U$ is
just $(Q\vec{n})\otimes (Q\vec(n))^T$, where $Q$ is the transition
matrix of the corresponding dynamics on complete network,
$\vec{n}=M(\vec{l})\vec{\bold{1}}$ is the corresponding
node-based macrostate. Besides, $WM(\vec{l})W^T$ restricted on $V$
has zero effect on the dynamics of $\vec{n}$. So $WM(\vec{l})W^T$
preserve the partial order of $\vec{n}$. This completes the proof of
theorem 2.


\section{Application}

The many significant consequences of monotone dynamical system,
especially low-dimensional ones have been explored in \cite{Hirsch,Smith,
Smith2}. In particular this substantial reduction in complexity of
the phase trajectories and the organization of the phase space into
hyperbolic equilibria and the heteroclinic orbits that connect them
have clear implications for mathematical sociology and network
science.

For application of monotonicity, one additional property \emph{\textbf{approximate from below (above)}} is important.\\
  \textbf{Definition:} $x$ can be approximated from below (above), if there is
a sequence $\{x_n\}$ satisfying $x_n\prec x_{n+1}\prec x$ ($x_n\succ
x_{n+1}\succ x$) for $n\geq 1$ and $x_n\go x$ as $n\go \infty$.

With the partial order for signalling system we discussed in this
paper, every point $x$ in macrostate space can be approximated from
below and above except for the two consensus state, since the
approximating sequence is given by $x+\epsilon_n\vec{e_i}$ with
$\vec{e_i}\succ 0$ and $\epsilon_n\go 0_- (0_+)$. Besides, the
macrostate space is a simplex therefore is finite-dimensional,
convex and compact.

Considering the properties of signalling system mentioned above, the
discussion in \cite{Smith} implies the following relevant
consequences:

\textbf{Theorem 3 - Global asymptotic stability \cite{Smith,Smith2}:} If a monotone signalling system contains exactly one equilibrium $e$, then every initial macrostate converge to $e$.

\textbf{Theorem 4 - Tipping point \cite{Smith,Smith2}:} If a monotone signalling
system has two equilibria $x\prec y$, $[x,y]$ denotes the set of all the points $z$ that $x\prec z \prec y$, then one of the following holds:\\
1) $y$ is stable, every point except $x$ converges to $y$.\\
2) $x$ is stable, every point except $y$ converges to $x$.\\
3) $x, y$ are both stable, and there exists another equilibrium $z\in [x,y]$, $z\neq x,y$.

According to the above two theorems, the global convergence of a
monotone signalling system is simply decided by the equilibria.
Theorem 3 guarantees the global convergence without knowing the
stability of the equilibrium. Theorem 4 is especially relevant to
the ``tipping point'' phenomenon found in 2-word Naming Game
\cite{ZhangLim},\cite{ZhangLimBolek12}, \cite{Xie11}. It predicts
all possible global structures of the 2-word NG dynamics from only
their monotonicity property instead of detailed inter-agent rules.
Therefore the previous results obtained for the NG
\cite{ZhangLim},\cite{ZhangLimBolek12}, \cite{Xie11},
\cite{Baronchelli12} can now be qualitatively generalized to any
monotone binary-signalling systems with two consensus states. If the
signalling system contains more than two ordered stable equilibria,
say $x\prec y \prec z$, then Theorem 4 can be applied on the domains
(called attracting basins) $[x, y]$ and $[y, z]$ separately.

\textbf{Theorem 5 - Convergence \cite{Smith,Smith2}}
For a monotone
signalling system, $M=\overline{Int(C)}$. Here $M$ is
the macrostate space, $C$ is the set of points that will eventually
converge. In another word, $Int(\overline{C^c})=\emptyset$, the set
of points that do not converge to anywhere is a nowhere dense set.

Considering that the signalling systems in the real world always
contain some noise, it is impossible for the trajectory to stay
inside a no-where dense set even with infinitely small noise,
therefore a monotone signalling system starting from any initial
state will eventually go to a stable equilibrium state.

Taken together, these theorems on Monotonicity and the new results
in this paper are applicable to many of the convergence, coherence
or synchrony questions that arise in mathematical biology and
ecology \cite{Levin}.


\section*{Acknowledgement}
This work was supported in part by the Army Research Office Grant
No. W911NF-09-1-0254 and W911NF-12-1- 467 0546. The views and
conclusions contained in this document are those of the authors and
should not be interpreted as representing the official policies,
either expressed or implied, of the Army Research Office or the U.S.
Government.

\newpage


\end{document}